\documentclass[12pt,reqno]{amsart}

\numberwithin{equation}{section}

\usepackage{nameref}

\usepackage{mdwlist}

\usepackage[final, 
color]{showkeys}
\definecolor{refkey}{gray}{.85}
\definecolor{labelkey}{gray}{.85}

\usepackage{AKstyle}

\usepackage[pagebackref=true, colorlinks]{hyperref}

\hypersetup{pdffitwindow=true,linkcolor=blue,citecolor=blue,urlcolor=blue,menucolor=blue}

\usepackage{comment}

\usepackage{caption}
\usepackage[labelfont=rm]{subcaption}

\makeatletter
\let\orgdescriptionlabel\descriptionlabel
\renewcommand*{\descriptionlabel}[1]{%
  \let\orglabel\label
  \let\label\@gobble
  \phantomsection
  \edef\@currentlabel{#1}%
  \let\label\orglabel
  \orgdescriptionlabel{#1}%
}
\makeatother

\newcommand\gdOvep{15}
\newcommand\gdO{17}
\newcommand\cCpow{10}
\newcommand\rhogb{\gb}
\newcommand\gbrho{\rho}

\begin{document}

\author{Jean Bourgain}
\thanks{Bourgain is partially supported by NSF grant DMS-0808042.}
\email{bourgain@ias.edu}
\address{IAS, Princeton, NJ}
\author{Alex Kontorovich}
\thanks{Kontorovich is partially supported by  
an NSF CAREER grant DMS-1254788, an Alfred P. Sloan Research Fellowship, and a Yale Junior Faculty Fellowship.}
\email{alex.kontorovich@yale.edu}
\address{Yale University, New Haven, CT}

\title[Thin Hypotenuses]
{The Affine Sieve Beyond Expansion I: Thin Hypotenuses}

\begin{abstract}
We study 
an 
instance of the Affine Sieve, 
 producing
a level of distribution 
beyond 
that which
can
be obtained from current techniques, 
even assuming a Selberg/Ramanujan-type spectral gap.
In particular, we consider
the set of 
hypotenuses 
in a
thin orbit 
of Pythagorean triples.
Previous work \cite{MyThesis, Kontorovich2009, KontorovichOh2012}
gave an exponent of distribution $\ga<1/12$ coming from Gamburd's  \cite{Gamburd2002} gap $\gt=5/6$, thereby producing $R=13$ almost primes in this linear sieve problem (see \S\ref{sec:intro} for definitions). 
If conditioned on a 
best possible gap
$\gt=1/2$, the known method 
would give
an exponent $\ga<1/4$, 
and
 $R=5$ almost primes.
The exponent $1/4$ is the natural analogue of the ``Bombieri-Vinogradov'' range of distribution for this problem, see Remark \ref{rmk:BVEH}.

In this paper, we
unconditionally
prove
 the exponent $\ga<7/24$ (in the ``Elliott-Halberstam'' range), 
thereby producing
$R=4$ almost primes. 
The main tools involve
 developing bilinear forms and the dispersion method in the range of incomplete sums for this Affine Sieve problem. 
\end{abstract}
\date{\today}
\maketitle
\tableofcontents

\newpage

\section{Introduction}\label{sec:intro}

\subsection{Level of Distribution
}\

The purpose of this paper is to inject bilinear forms, the dispersion method, and incomplete sums into the Affine Sieve to 
improve on known levels of
 of distribution 
 for a certain
 thin set of integers.
By this we mean the following.

We call a set $\cS\subset\Z^{n}$ 
 {\it thin} if there is some $\vep>0$ so that 
$$
\#(\cS\cap B_{N})<\#(\Zcl(\cS)\cap \Z^{n}\cap B_{N})^{1-\vep},
$$
for all $N$ large; here $B_{N}$ refers to the Euclidean ball of radius $N$ about the origin in $\R^{n}$, and $\Zcl(\cS)$ is the Zariski closure of $\cS$.
When $n=1$ and $\cS\subset\Z$ is infinite, the Zariski closure $\Zcl(\cS)=\bA$ is just affine space, so a set of integers being thin means that
$$
\#(\cS\cap[-N,N])
<N^{1-\vep}
,
$$
for $N$ large.
%

We say the set $\cS\subset\Z$ is {\it of Affine Sieve type} if there exists a triple $(\G,\bx_{0},f)$ with
\begin{itemize}
\item 
$\G\subset \GL_{n}(\Z)$ 
a
finitely generated,
not virtually
abelian
semigroup,
\item
$\bx_{0}\in\Z^{n}$ a primitive vector, giving rise to the orbit 
$$
\cO:=\bx_{0}\cdot \G
\subset \Z^{n}
,
$$
and
\item
$f:\Z^{n}\to\Z$ a polynomial,
\end{itemize}
so that 
$$
\cS=f(\cO).
$$

\begin{rmk}
 The orbit $\cO$ could be  thin
 without the set $\cS$ being 
 so.
For just two examples, the Apollonian orbit is thin, while the set of curvatures is not (see \cite{Kontorovich2013, BourgainKontorovich2012} for definitions and statements); likewise the 
orbit
 of absolutely Diophantine fractions is thin, while recent progress on Zaremba's conjecture \cite{BourgainKontorovich2011, BourgainKontorovich2011a} shows that the set of corresponding denominators is not.
\end{rmk}


Let $\cA=\{a_{N}(n)\}_{n\in\Z}$ be a sequence of nonnegative numbers supported on
$\cS
\cap[-N,N]
$.
%
Let $\cX$ be roughly the full ``mass'' of this sequence, that is, a quantity satisfying
$$
\cX\sim
|\cA|:=\sum_{n}a_{N}(n)
.
$$
For a square-free integer $\fq\ge1$, let
$$
|\cA_{\fq}|:=\sum_{n\equiv0(\fq)}a_{N}(n).
$$

If $\cS$ does not often favor multiples of some integers over others, then 
one might expect, for $\fq$ not too large relative to $N$, that 
\be\label{eq:cAq}
|\cA_{\fq}|=\rhogb(\fq)\cX+r(\fq)
.
\ee
Here the ``remainder'' $r(\fq)$ should be thought of as an error term, and $\rhogb(\fq)$ is a multiplicative ``local density,'' which 
in 
our
application will be roughly $1/\fq$ on average.\footnote{%
More generally, $\rhogb(p)$ could be about $\gk/p$ on average, for a fixed constant $\gk$ called the {\it sieve dimension}; our applications here will deal with a linear ($\gk=1$) sieve.}
In particular, we require that for any $2\le w<z$,
\be\label{eq:gbBnd}
\prod_{ w\le p<z}
(1-\rhogb(p))^{-1}
\le
C\cdot {\log z\over \log w}
,
\ee
for some $C>0$. 

Then a {\it level of distribution} 
is a number $\cQ$ so that the
total error $\cE$, that is,
the  remainders $|r(\fq)|$ summed up to $\cQ$, still does not exceed the full mass:
\be\label{eq:rqlcA}
\cE:=\sum_{\fq<\cQ}|r(\fq)|<\cX^{1-\vep}
.
\ee
Note that the level $\cQ$ is not intrinsic to the set $\cS$, but is instead a function of what one can prove about $\cS$. A key observation in the Affine Sieve, pioneered by 
Bourgain-Gamburd-Sarnak \cite{BourgainGamburdSarnak2006,BourgainGamburdSarnak2010}, is
that if $\G$ is an {\it expander} (that is, has a uniform spectral gap over congruence towers, see \eqref{eq:specGap}), then the sequence $\cS$ has a level of distribution 
\be\label{eq:QTalph}
\cQ=N^{\ga}, 
\ee
for some  $\ga>0$, called the  {\it exponent of distribution}.

We now describe our particular sequence $\cS$.

 \begin{figure}
        \begin{subfigure}[t]{3in}
                \centering
                {\includegraphics[height=1.5in]{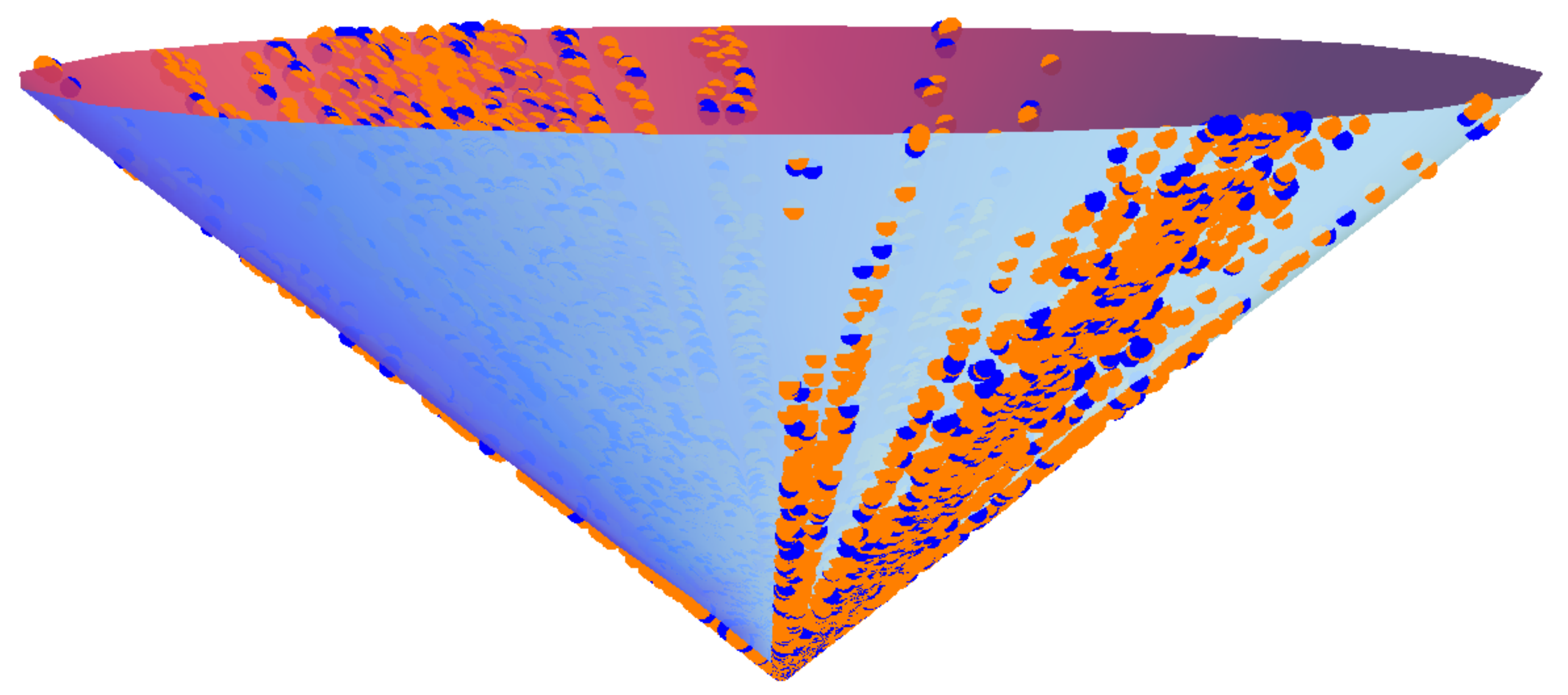}}
                \caption{View from the side}
                \label{fig:PythA}
        \end{subfigure}%
\qquad
        \begin{subfigure}[t]{1.5in}
                \centering
                \includegraphics[height=1.5in]{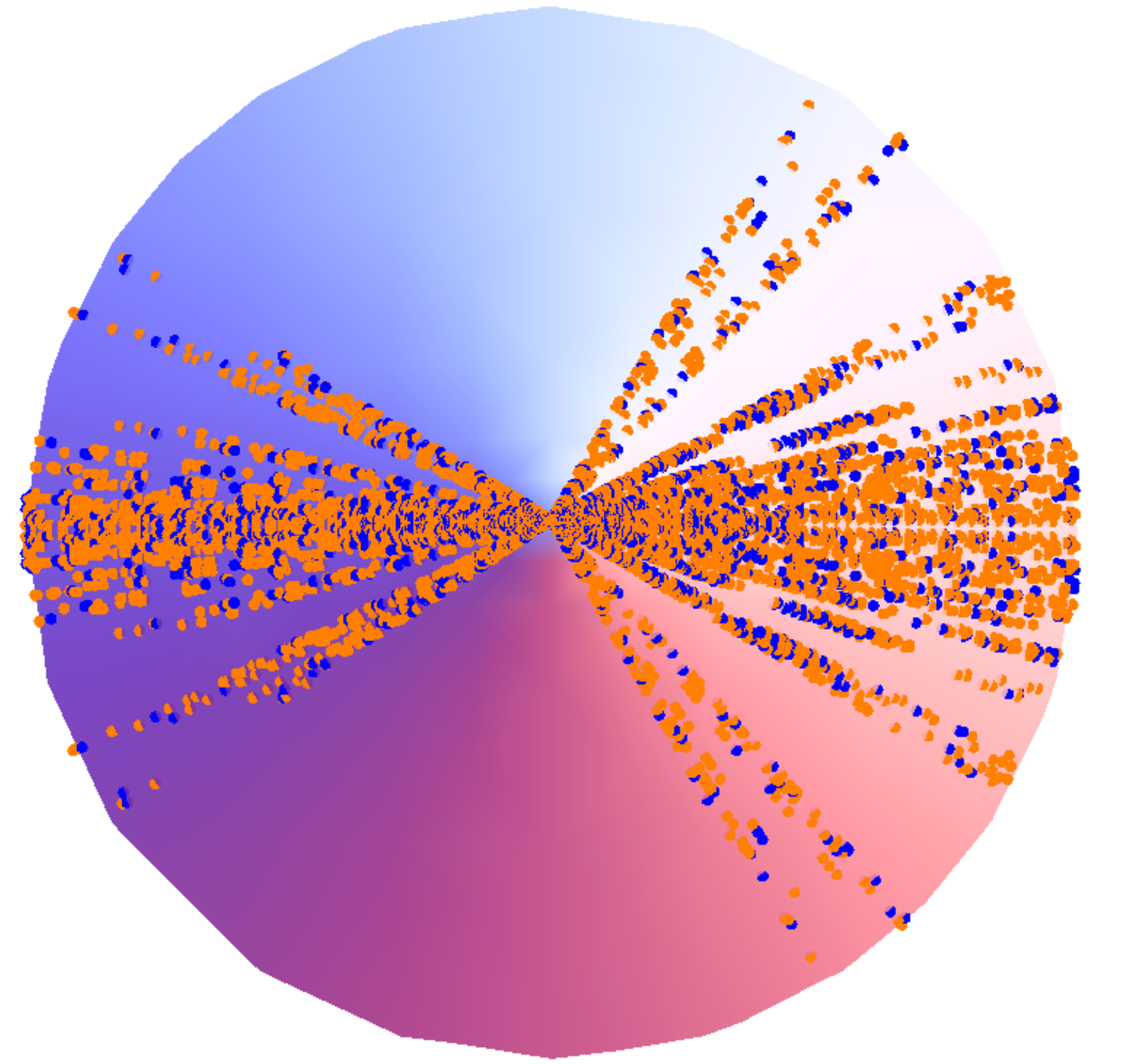}
                \caption{View from below}
                \label{fig:PythB}
        \end{subfigure}

\caption{A thin
Pythagorean 
orbit $\cO$. Points are marked according to whether the
 hypotenuse 
 is 
 prime (\protect\includegraphics[width=.1in]{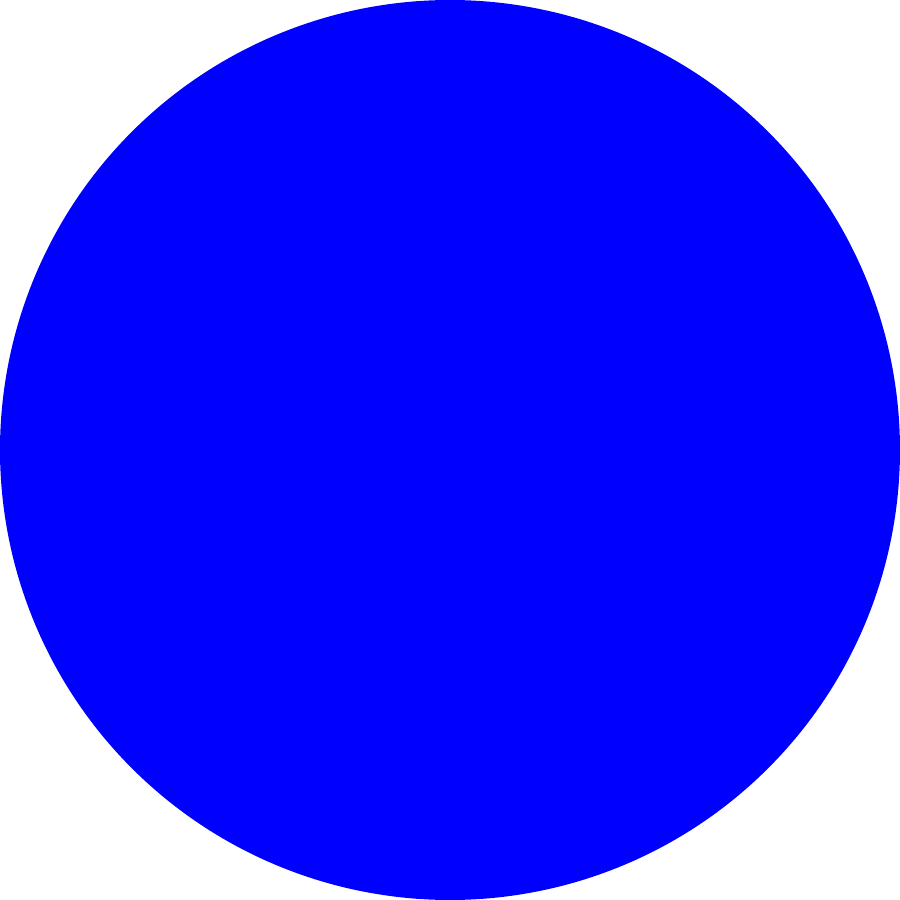}) or 
  composite (\protect\includegraphics[width=.1in]{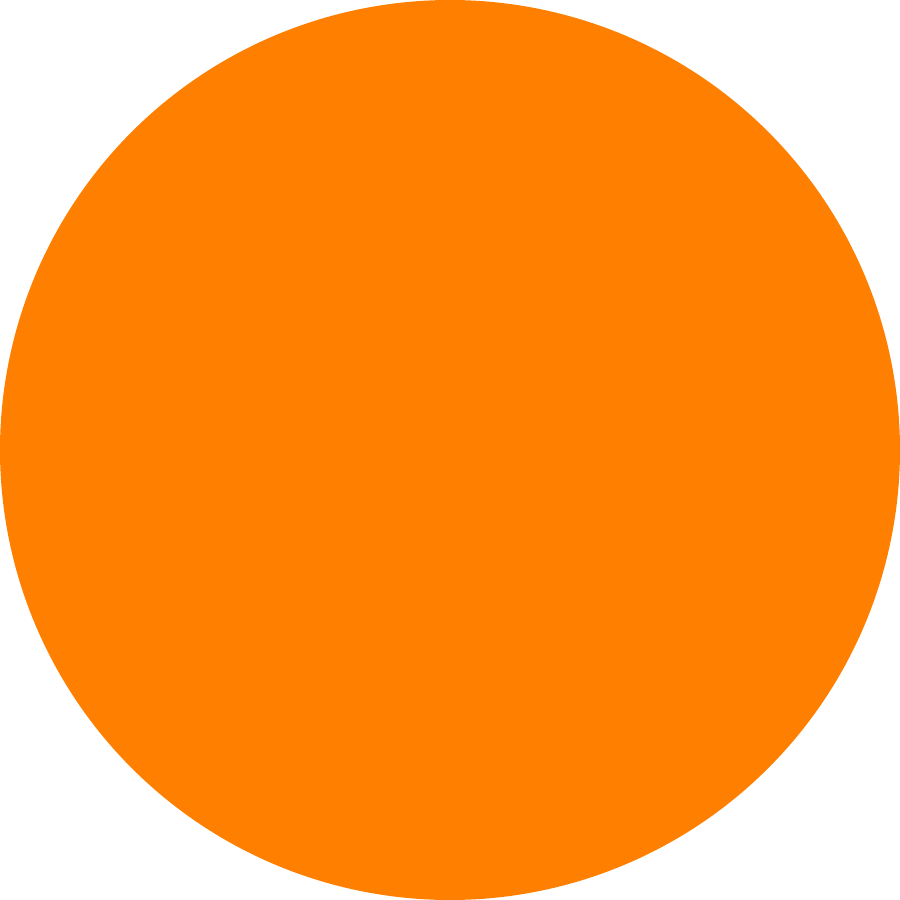})
  .}
\label{fig1}
\end{figure}

\subsection{Thin Orbits of Pythagorean Triples}\

A {\it Pythagorean triple} $\bx=(x,y,z)\in\Z^{3}$ is an integral point on the cone $F=0$, where 
$$
F(\bx):=x^{2}+y^{2}-z^{2}
.
$$ 
Let 
$$
G:=\SO^{\circ}_{F}(\R)=\SO^{\circ}(2,1)
$$ 
be the connected component of the identity of the special orthogonal group preserving $F$, and let 
$$
\G<\SO^{\circ}_{F}(\Z)
$$
be any 
geometrically finite 
subgroup  with integer entries. Assume further that $\G$ has no unipotent elements besides $I$ 
(otherwise, classical tools are available in the problem below). For a fixed Pythagorean triple $\bx_{0}$, say $\bx_{0}=(1,0,1)$, consider the orbit 
$$
\cO:=\bx_{0}\cdot\G, 
$$ 
and let $f(\bx)$ be the ``hypotenuse'' function,
$$
f(x,y,z):=z.
$$
Then we take the set $\cS$ of integers to be 
$$
\cS=f(\cO),
$$
that is, $\cS$ is the set of ``hypotenuses'' of triples in $\cO$. This set is clearly of Affine Sieve type. 
A sample such orbit is illustrated in Figure \ref{fig1}, where triples are marked according to whether or not their hypotenuses are prime.

Let 
$$
\gd=\gd_{\G}\in[0,1]
$$
be the critical exponent of $\G$, that is, the abscissa of convergence of the Poincar\'e series for $\G$; equivalently, $\gd$ is the Hausdorff dimension of the limit set of $\G$.
We will soon assume that $\gd$ is very near $1$, so consider henceforth the case $\gd>1/2$.
Then the set $\cS$ is thin if and only if $\gd<1$; indeed, it is 
known
 \cite{MyThesis, Kontorovich2009, KontorovichOh2012} 
that
$$
\#(\cS\cap[-N,N] )
\sim C\cdot 
N^{\gd}.
$$

\begin{thm}\label{thm:main}
Fix notation as above. Then for any $\vep>0$, there is some $\gd_{0}=\gd_{0}(\vep)<1$ with the following property. 
Whenever $\G$ has exponent $\gd>\gd_{0}$,
there exists a sequence $\cA=\{a_{N}(n)\}$ supported on $\cS$ 
so that \eqref{eq:cAq}--
\eqref{eq:rqlcA} hold, with exponent of distribution
\be\label{eq:gaIs}
\ga=7/24-\vep.
\ee
%
\end{thm}

See \S\ref{sec:cAis} 
for a precise construction of $\cA$. 

\begin{cor}\label{cor:main}
There is an absolute $\gd_{0}<1$ so that if $\gd>\gd_{0}$, then the set $\cS$ in Theorem \ref{thm:main} contains an infinitude of $R=4$ almost primes.
\end{cor}

\begin{rmk}\label{rmk:gaMore}
It may be possible to push the method to prove the exponent $\ga=7/24+\vep_{0}$ for a tiny $\vep_{0}>0$, see Remark \ref{rmk:vep0}. We have no applications for this improvement, so do not pursue it here.
\end{rmk}

\begin{rmk}\label{rmk:prev}
The exponent $\ga$ in \eqref{eq:gaIs} improves on the best previously available  exponent $\ga=1/12-\vep$, produced in \cite{MyThesis, Kontorovich2009, KontorovichOh2012}. Conditioned on an optimal spectral gap, the limit of that method gives $\ga=1/4-\vep$, see \S\ref{sec:method}. 
A key feature of our method is to divorce  the exponent of distribution from the spectral gap. Indeed,
were we to assume unproved hypotheses on infinite volume spectral gaps, the final value of our $\ga$ would not improve (though the value of $\gd_{0}$ would), see Remarks \ref{rmk:NeedGap} and \ref{rmk:NoSpecGap}.
\end{rmk}

\begin{rmk}\label{rmk:delVals}
We have made no effort to optimize the values of $\gd_{0}(\vep)$ and $\gd_{0}$ in Theorem \ref{thm:main} and  Corollary \ref{cor:main}, respectively. 
Our proof shows that 
$\gd_{0}(\vep)=1-10^{-\gdOvep}\vep$ suffices in Theorem \ref{thm:main}, and 
$\gd_{0}=1-10^{-\gdO}$ suffices for Corollary \ref{cor:main}.  These values can surely be improved.
\end{rmk}

Before explaining the source of 
this progress%
, we first reformulate and parametrize the problem. The group $\SO^{\circ}_{F}(\R)$ has a double cover by $\SL_{2}(\R)$, given explicitly by the map $\iota:\SL_{2}(\R)\to\SO^{\circ}_{F}:\mattwo abcd
{\mapsto}
$
$$
{1\over ad-bc}
\left(
\begin{array}{ccc}
 \frac{1}{2} \left(a^2-b^2-c^2+d^2\right) & c d-a b & \frac{1}{2}
   \left(-a^2-b^2+c^2+d^2\right) \\
 b d-a c & b c+a d & a c+b d \\
 \frac{1}{2} \left(-a^2+b^2-c^2+d^2\right) & a b+c d & \frac{1}{2}
   \left(a^2+b^2+c^2+d^2\right)
\end{array}
\right)
.
$$
Observe that, with $\bx_{0}=(1,0,1)$, we have
$$
f(\bx_{0}\cdot\iota\mattwo abcd)=c^{2}+d^{2}.
$$
Because the map $\iota$ is quadratic, a ball of radius $N$ in $\G$ is a ball of radius 
 $T$ in $\iota^{-1}(\G)$, where
\be\label{eq:T2N}
T^{2}\asymp N.
\ee
Abusing notation, we 
henceforth
call $G=\SL_{2}(\R)$, rename $\iota^{-1}(\G)$ to just $\G<\SL_{2}(\Z)$,
and let $\bx_{0}:=(0,1)$. Then we have the orbit $\cO=\bx_{0}\cdot\G$,
and
rename
\be\label{eq:fcdIs}
f(c,d)=c^{2}+d^{2},
\ee 
so that 
$$
\cS=f(\cO)
=
\left\{c^{2}+d^{2}:
\mattwo**cd\in\G
\right\}
.
$$ 
The set $\cS$ of integers is then the same as before. 
This
 reformulation
 is 
 just an easy 
 consequence
  of the ancient parametrization
$$
(x,y,z)=(c^{2}-d^{2},2cd,c^{2}+d^{2})
$$
of Pythagorean triples, in which the hypotenuse is a sum of two squares.

We now sketch the key new ideas which give the above claimed 
improvements.

\subsection{The Main 
Ingredients}\label{sec:method}\

This section is intended to be a heuristic discussion for the reader's convenience; statements are not made precisely. We first illustrate the ``standard'' Affine Sieve procedure used in \cite{MyThesis, Kontorovich2009,KontorovichOh2012}.

Recalling \eqref{eq:T2N}, we 
switch
to $T$ as our main parameter.
 A natural candidate for the sequence $\cA=\{a_{T}(n)\}$ is to take
\be\label{eq:aTnIs1}
a_{T}(n)=
\sum_{\g\in\G\atop\|\g\|<T}
\bo_{\{f(\bx_{0}\cdot\g)=n\}},
\ee
where $\|\cdot\|$ is the Frobenius norm, $\|\g\|^{2}=\tr(\g\,{}^{t}\g)$.
Then $a_{T}(n)$ 
is clearly 
supported on $\cS\cap[-N,N]$, with $N=T^{2}$,
and we have
\be\label{eq:cAqIs}
|\cA_{\fq}|=
\sum_{\g\in\G\atop\|\g\|<T}
\bo_{\{f(\bx_{0}\cdot\g)\equiv0(\fq)\}}
.
\ee

Let $\G(\fq)$ be the ``principal congruence'' subgroup of $\G$, that is, the kernel of the mod $\fq$ projection map. Then $\G(\fq)$ is still a thin group, but has finite index in $\G$. By Strong Approximation, we may assume the projection is onto, so $\G/\G(\fq)\cong\SL_{2}(\fq)$. Moreover, let $\G_{\bx_{0}}(\fq)$ be the stabilizer of $\bx_{0}$ mod $\fq$,
\be\label{eq:Gbx0q}
\G_{\bx_{0}}(\fq):=\{\g\in\G:\bx_{0}\cdot\g\equiv\bx_{0}(\fq)\}.
\ee
Clearly $\G(\fq)<\G_{\bx_{0}}(\fq)<\G$. Then \eqref{eq:cAqIs} can be decomposed as
\be\label{eq:cAqDecomp}
|\cA_{\fq}|
=
\sum_{\g_{0}\in\G_{\bx_{0}}(\fq)\bk\G}
\bo_{\{f(\bx_{0}\cdot\g_{0})\equiv0(\fq)\}}
\left[
\sum_{\g\in\G_{\bx_{0}}(\fq)}
\bo_{\|\g\g_{0}\|<T}
\right]
.
\ee
The inner sum (
suitably
smoothed
) is analyzed by spectral and representation-theoretic methods (see Theorem \ref{thm:specCount}), which prove modular equidistribution in essentially the following form:
$$
\Bigg[
\cdots
\Bigg]
=
\frac1{[\G:\G_{\bx_{0}}(\fq)]}\
C\cdot  T^{2\gd}
+
O(T^{2\gt
})
.
$$
Here $\gt<\gd$ is a spectral gap for $\G$ (see \S\ref{sec:spec}),
and both $\gt$ and the implied constant
are independent of $\fq$ and $\g_{0}$.
 For example, if $\gd>5/6$, then $\gt=5/6
$ 
is known
\cite{Gamburd2002}, whereas $\gt=1/2
$ would be a Selberg/Ramanujan quality gap (such a gap can be false in general). It is easy to compute the index $[\G:\G_{\bx_{0}}(\fq)]\sim 
\fq^{2}$, and the number of $\g_{0}$ with $f(\bx_{0}\cdot\g_{0})\equiv0(\fq)$ is about $\fq$. We thus have
\eqref{eq:cAq} and \eqref{eq:gbBnd}, with 
\be\label{eq:cXis}
\cX\asymp T^{2\gd}=N^{\gd},
\ee
and 
$$
|r(\fq)|\ll \fq T^{2\gt}.
$$
Then the level of distribution $\cQ$ is determined by requiring that
$$
\cE
=
\sum_{\fq<\cQ}|r(\fq)|
\ll
\cQ^{2}T^{2\gt
}
=
\cQ^{2}N^{\gt
}
$$
be an arbitrarily small power less than $\cX$. Compared to \eqref{eq:cXis}, we can take $\cQ=N^{\ga}$ with $\ga$ almost as large as $(\gd-\gt)/2$. Assuming $\gd$ is very near $1$ and applying Gamburd's gap $\gt=5/6$ gives an exponent $\ga$ almost as large as $1/12$. 
Under Selberg/Ramanujan, the biggest we could hope to make $\gd-\gt$ is just 
below
$1/2$, giving the conditional exponent $\ga<1/4$, as claimed in Remark \ref{rmk:prev}.

\begin{rmk}\label{rmk:BVEH}
Recall that for the sequence of primes in an arithmetic progression, the exponent of distribution $1/2$ follows from the Generalized Riemann Hypothesis, and the celebrated Bombieri-Vinogradov Theorem recovers this exponent unconditionally. The Elliott-Halberstam Conjecture predicts the exponent $1-\vep$, and some spectacular applications follow from {\it any} improvement on the exponent $1/2$. Returning to our sequence, the elements in $\G$ are taken of size $T$, and, as just explained, Selberg/Ramanujan would give the level $\cQ$ almost as large as $T^{\frac12}=N^{\frac14}$. So for this Affine Sieve problem, the exponent $1/4$ is the natural analogue of the ``Bombieri-Vinogradov'' range, and any exponent exceeding $1/4$ can be considered in the ``Elliott-Halberstam'' range. Theorem \ref{thm:main} produces just such an exponent, with $T^{\frac12}$ improved to almost  $T^{\frac12+\frac16}$. 
\end{rmk}

We now outline the new 
ingredients
introduced 
to prove Theorem \ref{thm:main}. Instead of the decomposition \eqref{eq:cAqDecomp}, we 
convert the problem into one on abelian harmonics, which are 
often
better understood. To this end, write
$$
\bo_{\{n\equiv0(\fq)\}}
=
\frac1\fq
\sum_{b(\fq)}
e_{\fq}(bn)
=
\frac1\fq
\sum_{q\mid\fq}
\sideset{}{'}\sum_{b(q)}
e_{q}(bn)
,
$$
where we have decomposed into primitive harmonics. Then 
\be\label{eq:cAqfq}
|\cA_{\fq}|
=
\sum_{\g\in\G\atop\|\g\|<T}
\frac1\fq
\sum_{q\mid\fq}
\sideset{}{'}\sum_{b(q)}
e_{q}(bf(\bx_{0}\cdot\g_{0}))
,
\ee
and we could 
try
 to 
break
the $q\mid\fq$ sum according to whether the modulus $q$ is below or above some parameter $Q_{0}<\cQ$. In the low range, we apply the standard spectral procedure as before. For $q$ large, we hope that there is sufficient cancellation already to treat the entire contribution as error. 
Summing these 
terms
up to the level $\cQ$, we need to estimate 
an exponential sum
essentially of the form
$$
\cE=
\sum_{Q_{0}<q<\cQ}\sum_{\g\in\G\atop\|\g\|<T}\frac1q\sideset{}{'}\sum_{b(q)} e_{q}(bf(\bx_{0}\cdot\g))
.
$$
The trivial bound here is $\cX\cQ$, so we need to save a tiny power more than $\cQ$. It seems hopeless to estimate the $\g$ sum over the 
intractable
thin group $\G$, so we use Cauchy-Schwarz to get rid of it. Recalling the parametrization \eqref{eq:fcdIs}, we estimate
\be\label{eq:cE2bnd}
\cE^{2}\ll \cX \cdot
\sum_{|c|,|d|< T}
\left|
\sum_{Q_{0}<q<\cQ}\frac1q\sideset{}{'}\sum_{b(q)} e_{q}(b(c^{2}+d^{2}))
\right|^{2}
.
\ee
It is here that we 
have made
critical use of the assumption that $\G$ has no parabolic elements: the bottom row $(c,d)$ of $\g$ determines $\g$ uniquely, 
allowing us to
 extend the $(c,d)$ sum to all of $\Z^{2}$. In so doing, we have essentially replaced the thin group $\G$ by all of $\SL_{2}(\Z)$, a loss we can only overcome if the dimension $\gd$ is 
 at least some $\gd_{0}$
 sufficiently close to $1$.
 
 Having squared
 , we now need to save a bit more than $\cQ^{2}$ from $\cE^{2}$. 
Unfortunately, the diagonal in \eqref{eq:cE2bnd}  contributes 
at most
a savings of $\cQ^{2}$ and no more, so this is hopeless. 
Taking a cue from Vinogradov's bilinear forms methods, as we did in \cite{BourgainKontorovich2010, BourgainKontorovich2011a, BourgainKontorovich2012}, we return to the sequence $\cA$, replacing \eqref{eq:aTnIs1} by
$$
a_{T}(n)=
\sum_{\g\in\G\atop\|\g\|<X}
\sum_{\gw\in\G\atop\|\gw\|<Y}
\bo_{\{f(\bx_{0}\cdot\g\gw)=n\}}
,
$$
where $XY=T$. Clearly $a_{T}$ is still supported on $\cS\cap[-N,N]$, but we are now taking greater advantage of the group structure of $\G$.  We treat $X$ as the long variable, being near $T^{1-\vep}$, and $Y$ small, of size $T^{\vep}$.
Following the above procedure, \eqref{eq:cE2bnd} is now replaced by
$$
\cE^{2}\ll X^{2\gd} \cdot
\sum_{|c|,|d|< T}
\left|
\sum_{\gw\in\G\atop\|\gw\|<Y}
\sum_{Q_{0}<q<\cQ}\frac1q\sideset{}{'}\sum_{b(q)} e_{q}(b f((c,d)\cdot\gw))
\right|^{2}
.
$$
Here
we have enough variables inside the square to hope to get cancellation. 
Expanding the square, the resulting modulus can be as large as $\cQ^{2}$ while the length of the $c,d$ sum is $T$, so to get beyond $\cQ=T^{1/2}=N^{1/4}$, we need to analyze incomplete sums. This is more-or-less standard, but 
unfortunately it is still impossible to get the desired cancellation. The issue now is that, on computing  the relevant 
Gauss sums, there is insufficient cancellation due to a lower order ``main'' term. This is a familiar feature in Linnik's 
work on
quadratic forms, and we must 
develop
a variant on his dispersion method here. 
Instead of just decomposing
\eqref{eq:cAqfq} according to $q<Q_{0}$ or $q\ge Q_{0}$, we 
add and subtract off this
lower order main term, see \S\ref{sec:Dispersion}. 

This change almost does the job but still not quite. We need to save in the $\gw$ sum (over the thin group $\G$) with a certain modular restriction for all moduli $q<\cQ$, while $\gw$ is only of small size $Y$. To do so, we make one 
final
technical modification, breaking the $\gw$ sum into two, one of 
size $Y_{1}$ and another of much smaller size $Y_{2}$, 
with $Y_{1}Y_{2}=Y$. 
For $q$ not too large, we apply spectral theory in $Y_{1}$. For larger $q$, we use the fact that $Y_{2}$ is very much smaller to turn the modular restriction into an archimedean
one,
 giving the desired savings; see Theorem \ref{thm:gWsave}.






\subsection{Outline}\

The paper proceeds as follows. We use \S\ref{sec:AutF} to recall some facts about infinite volume spectral and representation theory, spectral gaps, and decay of matrix coefficients. These are used to prove certain counting statements needed in the sequel. In \S\ref{sec:setup}, we 
construct
the sequence $\cA$ and divide it into a main term and error. The former is analyzed in \S\ref{sec:main} and the latter in \S\ref{sec:error}. These estimates are collected in \S\ref{sec:proofs} to complete the proofs of Theorem \ref{thm:main} and Corollary \ref{cor:main}.

\subsection{Notation}\

We use the following standard notation. Define $e(x):=e^{2\pi i x}$ and set $e_{q}(x):=e(\frac xq)$. We use $f\ll g$ and $f=O(g)$ interchangeably; moreover $f\asymp g$ means $f\ll g\ll f$.
Unless otherwise specified, the implied constants may depend at most on $\G$
, which is treated as fixed.
The letter $\vep$ always denotes an arbitrarily small positive constant; when it appears in an equation, the implied constant may also implicitly depend on $\vep$.
The letter $C$ is a 
positive constant, not necessarily the same at each occurrence. 
The letter $p$ always represents a prime number.
 The symbol $\bo_{\{\cdot\}}$ is the indicator function of the event $\{\cdot\}$. 
The greatest common divisor of $n$ and $m$ is written $(n,m)$ and their least common multiple is $[n,m]$
. The cardinality of a finite set $S$ is denoted $|S|$ or $\# S$.
The transpose of a matrix $g$ is written ${}^{t}g$. The prime symbol $'$ in $\underset{b(q)}{\gS} {}'$ means the range of $b(\mod q)$ is restricted to $(b,q)=1$. 

\subsection*{Acknowledgements} The authors are grateful to Henryk Iwaniec and Peter Sarnak for illuminating discussions.

\newpage

\section{Background: Infinite Volume Automorphic Forms}\label{sec:AutF}

\subsection{Spectral Theory, Expansion, and Matrix Coefficients}\label{sec:spec}\

Let $\G<\SL_{2}(\Z)$ be a finitely generated thin subgroup with 
critical 
exponent
$\gd>1/2$. The hyperbolic Laplacian $\gD=-y^{2}(\dd_{xx}+\dd_{yy})$ acts on $L^{2}(\G\bk\bH)$, where $\bH=\{x+iy:x\in\R,y>0\}$ is the hyperbolic upper half plane. 
Then the spectrum $\Spec(\G)$ is purely continuous above $1/4$  and  consists of a finite number of discrete eigenvalues below $1/4$ \cite{LaxPhillips1982}. Labeling the discrete eigenvalues as
$$
0<\gl_{0}<\gl_{1}\le\cdots\le\gl_{max}<1/4,
$$
we have the Patterson-Sullivan formula for the base eigenvalue $\gl_{0}=\gd(1-\gd)$ \cite{Patterson1976, Sullivan1984}. 

For a square-free integer $q\ge1$, let $\G(q)$ be the principal congruence subgroup of $\G$. Write the discrete spectrum of  $\G(q)$ as
\be\label{eq:gls}
0<\gl_{0}(q)<\gl_{1}(q)\le\cdots\le\gl_{max(q)}(q)<1/4.
\ee
The inclusion $\G(q)<\G$ induces the reverse inclusion
$$
\Spec(\G(q))\supset\Spec(\G),
$$ 
and we call $\Spec^{new}(q)$ the complementary set of ``new'' spectrum of level $q$.
It is easy to see that
 $\G(q)$ 
 and
 $\G$ 
  have 
   the same base eigenvalue $\gl_{0}(q)=\gl_{0}$,
but the 
next
eigenvalue $\gl_{1}(q)$ could {\it a priori} approach $\gl_{0}$ as $q$ increases.
We say $\G$ is an {\it expander} if this doesn't happen. 
More precisely, we say $\G$ has a uniform {\it spectral gap} 
$$
 \gt\in(1/2,\gd),
$$
if
 there exists a number 
 \be\label{eq:fBspec}
 \fB\ge1
 \ee 
 so that for all $(q,\fB)=1$, the new spectrum lies above $\gt(1-\gt)$, 
\be\label{eq:specGap}
 \Spec^{new}(q)\subset[\gt(1-\gt),\infty)
 .
 \ee 
\begin{thm}[\cite{BourgainGamburdSarnak2010, BourgainGamburdSarnak2011, Gamburd2002}]\label{thm:spec}
Every $\G$ as above has some spectral gap $\gt<\gd$. Moreover if $\gd>5/6$, then $\G$ has the 
absolute
spectral gap $\gt=5/6$.
\end{thm}

The group $G=\SL_{2}(\R)$ acts by the right regular representation on $V:=L^{2}(\G(q)\bk G)$. 
By the Duality Theorem \cite{GelfandGraevPS1966},
we have the decomposition corresponding to \eqref{eq:gls} and \eqref{eq:specGap}:
\be\label{eq:Vdecomp}
V
=
\bigoplus
_{\gl_{j}<\gt(1-\gt)}
V_{\gl_{j}}
\quad
\oplus\quad
V^{\perp}
.
\ee
Here
each $V_{\gl_{j}}$ is a complementary series representation of parameter 
$1/2<s_{j}<1$, where $\gl_{j}=s_{j}(1-s_{j})$,
and
 $V^{\perp}$ does not weakly contain any complementary series representation of parameter $s>\gt$.

The following theorem on the decay of matrix coefficients is well-known \cite{CowlingHaagerupHowe1988}.
\begin{thm}
 For any right $K$-invariant vectors $\Psi_{1},\Psi_{2}\in V^{\perp}$,
\be\label{eq:decayMtrx}
|\<\pi(g).\Psi_{1},\Psi_{2}\>|
\ll
\|g\|^{2(1-\gt)}
\|\Psi_{1}\|\,\|\Psi_{2}\|
,
\ee
as $\|g\|\to\infty$.
\end{thm}

\subsection{Counting Statements}\

We give here some infinite volume counting statements needed in the sequel. These are more-or-less standard, but we give sketches for the reader's convenience.
Our first goal is to show that, in smooth form, we can count norm balls in $\G$, uniformly in congruence towers with sharp rates.
Once and for all, fix a smooth bump function $\psi$ on $G/K$, 
that is, we assume
 $\int_{G/K}\psi=1$,
 $\psi$ is non-negative,  and supported in a ball of radius $1/100$, say, about the origin. Then for $X>1$ and $g\in G$,
 the integral
\be\label{eq:gUXdef}
\gU_{X}(g):=
\int_{G/K}
\int_{G/K}
\bo_{\|h_{1}^{-1}gh_{2}\|<X}
\psi(h_{1})
\psi(h_{2})
dh_{1}dh_{2}
\ee
is well-defined, since our norm is bi-$K$-invariant.
It is easy to see that
$$
\gU_{X}(g)=
\threecase
{1}{if $\|g\|<\frac9{10}X$,}
{0}{if $\|g\|>\frac{11}{10}X$,}
{\in[0,1]}{otherwise,}
$$
so $\gU_{X}$ 
is a smoothed version of the indicator function $\bo_{\|g\|<X}$. 

\begin{thm}\label{thm:specCount}
Let $\G$ have exponent $\gd>1/2$ and spectral gap $\gt<\gd$, as above.
Then we have
\be\label{eq:normX}
\sum_{\g\in\G}
\gU_{X}(\g)
\sim C\cdot X^{2\gd},
\ee
as $X\to\infty$.
Moreover, for any $\g_{0}\in\G$,  any 
square-free $q$ coprime to $\fB$,
and any $\tilde\G(q)$ satisfying $\G(q)<\tilde\G(q)<\G$, we have
\be\label{eq:normXq}
\sum_{\g\in\tilde\G(q)}
\gU_{X}(\g\g_{0})
=
\frac1{[\G:\tilde\G(q)]}
\sum_{\g\in\G}
\gU_{X}(\g)
+
O(X^{2\gt})
.
\ee
The implied constant does not depend on $q$ or $\g_{0}$.
\end{thm}
\pf[Sketch of proof]
Let 
$$
\cF(g,h):=\sum_{\g\in\tilde\G(q)}\bo_{\{\|g^{-1}\g h\|<X\}}
,
\qquad
\Psi(g):=\sum_{\g\in\tilde\G(q)}\psi(\g g),
$$
and
$$
\Psi_{\g_{0}}(g):=\sum_{\g\in\tilde\G(q)}\psi(\g_{0}^{-1}\g g).
$$
Then $\Psi,\Psi_{\g_{0}}\in L^{2}(\tilde\G(q)\bk G/K)$, and $\cF\in L^{2}(\tilde\G(q)\bk G/K\times\tilde\G(q)\bk G/K)$.
After some changes of variables, unfolding and refolding integrals (see \cite[Lemma 3.7]{BourgainKontorovichSarnak2010}), we obtain
\be\label{eq:gUsum}
\sum_{\g\in\tilde\G(q)}
\gU_{X}(\g\g_{0})
=
\<\cF,\Psi\otimes\Psi_{\g_{0}}\>
=
\int_{K\bk G/K}
\bo_{\{\|g \|<X\}}
\<\pi(g).\Psi,\Psi_{\g_{0}}\>
dg
.
\ee
Expanding  spectrally according to 
\eqref{eq:Vdecomp}, we have
$$
\Psi=\sum_{\gl_{j}<\gt(1-\gt)}\<\Psi,\vf_{j}\>\vf_{j}
+
\Psi^{\perp}
,
$$
where $\vf_{j}$ is an $L^{2}$-normalized eigenfunction corresponding to $\gl_{j}$. Similarly expand $\Psi_{\g_{0}}$, and insert these expansions into the last inner product of 
\eqref{eq:gUsum}. 
Setting $q=1, \g_{0}=I$ and applying \eqref{eq:decayMtrx} gives \eqref{eq:normX} after a standard calculation. (Here we used  that $\|\Psi^{\perp}\|\le\|\Psi\|\ll1$, since the support of $\psi$ is absolute). For \eqref{eq:normXq}, we observe that $\vf_{j}$ are ``oldforms'', and that their normalization in $L^{2}(\tilde\G(q)\bk \bH)$ differs from that in $L^{2}(\G\bk\bH)$ by the factor $[\G:\tilde\G(q)]^{-1/2}$, whence the claim follows. 
\epf

Next we need to save a small power of $q$ for a modular restriction in a  ball, where $q$ can be much larger than the size of the ball. We accomplish this by replacing the ball by a product of two balls with vastly different sizes (a related  trick was used already in \cite[\S5]{BourgainKontorovich2012}).

Let $\G$ have exponent $\gd$ and spectral gap $\gt$ as above, and set
\be\label{eq:cCis}
\cC:={10^{\cCpow}\over \gd-\gt}.
\ee
Let $Y=Y_{1}Y_{2}$ with 
\be\label{eq:Y2toY1}
Y_{1}=Y_{2}^{\cC},
\ee
and let $\gW_{Y}$ denote the multi-subset of $\G$ given by
\be\label{eq:gWYis}
\gW_{Y}:=
\{
\gw_{1}\cdot\gw_{2}:\gw_{1},\gw_{2}\in\G,\|\gw_{1}\|<Y_{1},\|\gw_{2}\|<Y_{2}
\}
,
\ee
so that
\be\label{eq:gWYsize}
|\gW_{Y}|\sim 
C\cdot Y^{2\gd},
\ee
by \cite{LaxPhillips1982}.
Recall from \eqref{eq:fcdIs} that $f(c,d)=c^{2}+d^{2}$.
\begin{thm}\label{thm:gWsave}
Given
any sufficiently small $\eta>0$, 
 and any parameters
  $\cQ,Y\to\infty$, 
with
\be\label{eq:YtoT}
Y\ge\cQ^{
\eta\cdot\frac32\cdot{\cC+1\over \gd-\gt}
},
\ee 
 we have
the following.
For 
any square-free $q<\cQ$ with $(q,\fB)=1$, and any $\bx\in\Z^{2}$ with $
(\bx,q)=1
$, we have
\be\label{eq:gWYsave}
\sum_{\gw\in\gW_{Y}}
\bo_{\{f(\bx\cdot\gw)\equiv0(q)\}}
\ll q^{-\eta}|\gW_{Y}|
.
\ee
\end{thm}
\pf[Sketch of proof]
The proof decomposes into two cases, depending on the size of $q$.

{\bf Case 1: $q<Y_{1}^{\gd-\gt}$.}
Starting with the left side of \eqref{eq:gWYsave}, fix $\gw_{2}$, and decompose $\gw_{1}$ as
\beann
LHS\eqref{eq:gWYsave}
&=&
\sum_{\gw_{1}\in\G\atop\|\gw_{1}\|<Y_{1}}
\sum_{\gw_{2}\in\G\atop\|\gw_{2}\|<Y_{2}}
\bo_{\{f(\bx\cdot\gw_{1}\gw_{2})\equiv0(q)\}}
\\
&=&
\sum_{\gw_{2}\in\G\atop\|\gw_{2}\|<Y_{2}}
\sum_{\gw_{0}\in\G_{\bx}(q)\bk\G}
\bo_{\{f(\bx\cdot\gw_{0}\gw_{2})\equiv0(q)\}}
\left[
\sum_{\gw_{1}\in\G_{\bx}(q)}
\bo_{\{\|\gw_{1}\gw_{0}\|<Y_{1}\}}
\right]
,
\eeann
where $\G_{\bx}(q)$ is the subgroup of $\g\in\G$ for which $\bx\cdot\g\equiv\bx(q)$.
Applying Theorem \ref{thm:specCount} to the innermost sum and estimating, we get
$$
LHS\eqref{eq:gWYsave}
\ll
Y_{2}^{2\gd}
q
\left[
\frac1{q^{2}}
Y_{1}^{2\gd}
+
Y_{1}^{2\gt}
\right]
\ll
\frac1q
|\gW_{Y}|
,
$$
since $q<Y_{1}^{\gd-\gt}$. Thus in this range we can prove \eqref{eq:gWYsave} with $\eta=1$. 

{\bf Case 2: $q\ge Y_{1}^{\gd-\gt}$.} Using \eqref{eq:Y2toY1} and \eqref{eq:cCis}, we have in this range that
\be\label{eq:626ss}
q\ge Y_{1}^{\gd-\gt}=Y_{2}^{10^{\cCpow}}
.
\ee
Now we fix $\gw_{1}$ (with $\ll Y_{1}^{2\gd}$ choices) and play with $\gw_{2}$. We wish to use the discrepancy in the huge modulus $q$ relative to the small size $Y_{2}$ to convert the modular restriction into an archimedean one, as follows.

Set $\by=\bx\cdot\gw_{1}$, and note that $(\by,q)=1$, since $\det\gw_{1}=1$. 
Drop the subscript from $\gw_{2}$, and write $\by=(u,v)$, $\gw=\mattwos abcd$. We may assume without loss of generality that $(u,q)=1$, so that the condition $f(\by\cdot\gw)\equiv0(q)$ is equivalent to
\be\label{eq:626s}
(a^{2}+b^{2})+2(\bar u v)(ac+bd)+(\bar u v)^{2}(c^{2}+d^{2})\equiv0(q).
\ee
We need to estimate the cardinality of 
$$
\cT=\cT(Y_{2};u,v)
:=
\{
\gw\in\G:
\|\gw\|<Y_{2}
\text{ and \eqref{eq:626s} holds}
\}
.
$$
For each $\gw\in\cT$, let $P_{\gw}\in\Z[U,V]$ be the (linear) polynomial
$$
P_{\gw}(U,V):=a^{2}+b^{2}+U(ac+bd)+V(c^{2}+d^{2})
,
$$
and note that it has logarithmic height at most 
$$
h:=2\log Y_{2}.
$$ 
Consider the affine variety
$$
\cV:=\bigcap_{\gw\in\cT}
\{
P_{\gw}
=0\}
.
$$

We claim that $\cV(\C)\neq\O$, and argue by contradiction. If $\cV(\C)$ is empty, then effective Hilbert's Nullstellensatz \cite[Theorem IV]{MasserWustholz1983} gives the existence of polynomials $Q_{\gw}\in\Z[U,V]$ and an integer $\fd\ge1$ so that the Bezout equation 
\be\label{eq:626s3}
\sum_{\gw\in\cT}
P_{\gw}
\cdot
Q_{\gw}
=
\fd
\ee
holds, 
with $\fd\le \exp\left(8^{4\cdot 2-1}(h+8\log 8)\right)\ll Y_{2}^{10^{7}}$. 
(Better estimates exist, but this will suffice for our purposes.) In particular, 
 \eqref{eq:626ss} forces
\be\label{eq:fdq}
1\le\fd<q 
.
\ee
But reducing \eqref{eq:626s3} mod $q$ and setting $(U,V)\equiv(2\bar u v,(\bar uv)^{2})$ gives $\fd\equiv0(q)$ by \eqref{eq:626s}.
This   is incompatible with \eqref{eq:fdq}, giving our desired contradiction. Hence $\cV(\C)$ must be nonempty.

It is  easy to see  that the set of rational points $\cV(\Q)$ is then also non-empty, and hence, after clearing denominators, there exist coprime integers $t_{*},u_{*},v_{*}$ so that
$$
f_{*}(\gw):=
t_{*}(a^{2}+b^{2})+u_{*}(ac+bd)+v_{*}(c^{2}+d^{2}) 
=0
,
$$
for all $\gw\in\cT$. Hence we have finally lifted the modular restriction to an archimedean  one.

Now choose a prime 
$$
\ell\asymp Y_{2}^{2(\gd-\gt)/3},
$$ 
 replace $f_{*}(\gw)=0$ by the weaker condition $f_{*}(\gw)\equiv0(\ell)$, and proceed as before:
\beann
|\cT|
&\le&
\sum_{\gw\in\G\atop\|\gw\|<Y_{2}}
\bo_{\{f_{*}(\gw)\equiv0(\ell)\}}
=
\sum_{\gw_{0}\in\G(\ell)\bk\G}
\bo_{\{f_{*}(\gw_{0})\equiv0(\ell)\}}
\left[
\sum_{\gw\in\G(\ell)}
\bo_{\{\|\gw\gw_{0}\|<Y_{2}\}}
\right]
\\
&\ll&
\ell^{2}
\left[
\frac1{\ell^{3}}
Y_{2}^{2\gd}
+
Y_{2}^{2\gt}
\right]
\ll
\frac1\ell
Y^{2\gd}
,
\eeann
by assumption on the size of $\ell$. By \eqref{eq:Y2toY1} and \eqref{eq:YtoT}, we have thus saved
$$
\ell
\gg
Y_{2}^{2(\gd-\gt)/3}
=
Y^{\frac23\cdot{\gd-\gt\over \cC+1}}
\ge
\cQ^{\eta}
>
q^{\eta},
$$
as desired.
\epf

\newpage

\section{Setup, Construction of $\cA$, and Dispersion}\label{sec:setup}

\subsection{Initial manipulations}\

Recall that $\G$ is a thin, finitely generated subgroup of $\SL_{2}(\Z)$ with no parabolic elements and dimension $\gd>1/2$. Our sieve problem concerns the set $\cS$, where $\bx_{0}=(0,1)$, $\cO=\bx_{0}\cdot\G$, $f(c,d)=c^{2}+d^{2}$, and $\cS=f(\cO)$.
We first perform some initial manipulations
. 

By Strong Approximation, $\G(\mod p)$ is all of $\SL_{2}(p)$ except for a finite list $\cP$ of ``bad'' primes. 
We may increase $\cP$ if necessary to make sure that $2\in\cP$, and also that $\cP$ contains all the primes dividing $\fB$ in \eqref{eq:fBspec}.
Then renaming $\fB:=\prod_{p\in\cP}p$, it follows from 
 Goursat's Lemma that if $q$ is square-free with $(q,\fB)=1$, then 
 \be\label{eq:GmodQ}
 \G(\mod q)\cong\SL_{2}(q).
 \ee
At the cost of decreasing $\cS$, we may replace $\G$ by its principal congruence group $\G(\fB)$ of level $\fB$, renaming $\G$ and the resulting set $\cS$. Then observe that if $\cA=\{a_{N}(n)\}$ is supported on $\cS$ and $(\fq,\fB)>1$, then $|\cA_{\fq}|=0$, since there will be no $\mattwos abcd\in\G$ with $c^{2}+d^{2}\equiv0(p)$ for any $p\mid\fB$.
Similarly, $|\cA_{\fq}|$ vanishes if $\fq$ contains a prime factor $p\equiv3(\mod 4)$. We thus assume henceforth that any numbers $\fq, q$ are square-free, with prime divisors $p\equiv1(\mod 4)$ and $p\nmid\fB$.

Next we manipulate the function $f$.
To each $\g\in\G$, 
we attach the
 binary quadratic form 
$$
\ff_{\g}(x,y)=Ax^{2}+2Bxy+Cy^{2},
$$
where $\g\cdot{}^{t}\g=\mattwos ABBC$. Observe that if $\g=\mattwos **cd$, then $\ff_{\g}(0,1)=c^{2}+d^{2}$, so that 
$$
f(\bx_{0}\cdot\g)=\ff_{\g}(0,1).
$$
Moreover,  for another $\gw\in\G$,  we have
$$
\ff_{\g\gw}(0,1)=\ff_{\gw}(c,d).
$$

\subsection{Construction of $\cA$}\label{sec:cAis}\

We are now in position to construct our sequence $\cA$. 
Let $T$ be our main growing parameter, and 
write 
\be\label{eq:TtoXY}
T=XY,
\ee 
with parameters $X$ and $Y$ to be chosen in \S\ref{sec:proofs}. 
Recalling the smoothing function $\gU_{X}$ in \eqref{eq:gUXdef}, and the multi-set $\gW_{Y}$ in \eqref{eq:gWYis}, we
define 
$$
\cA=\{a_{T}(n)\}
$$ 
by
\bea
\nonumber
a_{T}(n):
&=&
\sum_{\g\in\G}
\gU_{X}(\g)
\sum_{\gw\in\gW_{Y}}
\bo_{\{f(\bx_{0}\cdot\g\gw)=n\}}
\\
\label{eq:aTDef}
&=&
\sum_{\g\in\G\atop\g=
\bigl( \begin{smallmatrix}
*&*\\ c&d
\end{smallmatrix} \bigr)
}
\gU_{X}(\g)
\sum_{\gw\in\gW_{Y}}
\bo_{\{\ff_{\gw}(c,d)=n\}}
.
\eea
We emphasize 
that $\gW_{Y}$ is a multi-set, so the $\gw$ sum in \eqref{eq:aTDef} is with multiplicity.

Since $f$ is quadratic, the support of $a_{T}(n)$ is in $\cS\cap[1,N]$ with 
\be\label{eq:NtoT}
N\asymp T^{2}.
\ee
By \eqref{eq:normX} and \eqref{eq:gWYsize}, we have crudely that
\be\label{eq:cAis}
|\cA|:=\sum_{n}a_{T}(n)\asymp T^{2\gd}.
\ee
For $\fq<\cQ<X<T$, set 
\be\label{eq:cAfqDef}
|\cA_{\fq}|:=\sum_{n}a_{T}(n)\bo_{\{n\equiv0(\fq)\}}.
\ee
Recall that $\fq$ is square free and
 a product of
 primes $p\equiv1(4)$, $p\nmid\fB$. 

\subsection{Setting Up the Dispersion Method}\label{sec:Dispersion}\

We apply a novel version of the dispersion method. Let $\gbrho(q)$ and $\Xi(q)$ be multiplicative functions, defined at $p$ by
\be\label{eq:gbDef}
\gbrho(p):={2p-1\over p^{2}}
,
\ee
and
\be\label{eq:XiDef}
\Xi(p;n):=\bo_{\{n\equiv0(p)\}}-\gbrho(p).
\ee
Here $\gbrho$ is a ``main term'' being subtracted off at each prime factor to make Lemma \ref{lem:Dispersion} hold.

Then inserting 
$$
\bo_{\{n\equiv0(\fq)\}}
=
\prod_{p\mid\fq}
\bo_{\{n\equiv0(p)\}}
=
\prod_{p\mid\fq}
\left(
\Xi(p;n)+\gbrho(p)
\right)
=
\sum_{q\mid\fq}
\Xi(q;n)\gbrho\left(\frac \fq q\right)
$$
into \eqref{eq:cAfqDef} gives
$$
|\cA_{\fq}|
=
\sum_{q\mid\fq}
\sum_{n}a_{T}(n)
\Xi(q;n)\gbrho\left(\frac \fq q\right)
.
$$

For a parameter $Q_{0}<\cQ<T$ to be chosen 
in \S\ref{sec:proofs},  we write 
\be\label{eq:cAqDec}
|\cA_{\fq}|=\cM_{\fq}+r(\fq)
,
\ee
say,
where we decomposed
according to whether $q<Q_{0}$ or $q\ge Q_{0}$. The ``main term'' 
\be\label{eq:cMfqIs}
\cM_{\fq}
:=
\sum_{q\mid\fq\atop q<Q_{0}}
\sum_{n}a_{T}(n)
\Xi(q;n)\gbrho\left(\frac \fq q\right)
\ee 
will be analyzed by spectral methods in the next section. 
Thereafter, we must control the net error 
\be\label{eq:cEDef}
\cE:=\sum_{\fq<\cQ}|r(\fq)|
\ee
up to level $\cQ$, with $\cQ$ as large as possible.

\newpage

\section{Analysis of $\cM_{\fq}$}\label{sec:main}

Keeping the previous notation, the goal of this section is to prove the following
\begin{thm}
Let $\gb$ be a multiplicative function defined on primes by
\be\label{eq:gbIs}
\gb(p)=
\twocase{}
{\frac2{p+1}}{if $p\equiv1(4)$, $p\nmid\fB$,}
{0}{otherwise,}
\ee
and let 
\be\label{eq:cXisis}
\cX:=|\cA|=|\gW_{Y}|\sum_{\g\in\G}\gU_{X}(\g).
\ee
Then there is a decomposition
\be\label{eq:cMfqD}
\cM_{\fq}=\gb(\fq)\cX+r^{(1)}(\fq)+r^{(2)}(\fq)
,
\ee
with
\be\label{eq:r1Bnd}
\sum_{\fq<\cQ}|r^{(1)}(\fq)|
\ll
\cQ^{\vep}
\cX
{
Q_{0}^{3}
\over
X^{2(\gd-\gt)}
}
.
\ee
and
\be\label{eq:r2Bnd}
\sum_{\fq<\cQ}|r^{(2)}(\fq)|
\ll
\cQ^{\vep}
\cX
\frac1{Q_{0}}
.
\ee
\end{thm}

\pf
Inserting \eqref{eq:aTDef} into \eqref{eq:cMfqIs}, we have that
$$
\cM_{\fq}
=
\sum_{q\mid\fq\atop q<Q_{0}}
\gbrho\left(\frac \fq q\right)
\sum_{\g\in\G}
\gU_{X}(\g)
\sum_{\gw\in\gW_{Y}}
\Xi(q;f(\bx_{0}\cdot\g\gw))
.
$$
Of course $\Xi(q;n)$ only depends on the residue class of $n(\mod q)$. Let $\G_{\bx_{0}}(q)$ be the stabilizer of $\bx_{0}(\mod q)$ in $\G$, as in \eqref{eq:Gbx0q}, and  decompose the $\g$ sum as
$$
\cM_{\fq}
=
\sum_{q\mid\fq\atop q<Q_{0}}
\gbrho\left(\frac \fq q\right)
\sum_{\gw\in\gW_{Y}}
\sum_{\g_{0}\in\G_{\bx_{0}}(q)\bk\G}
\Xi(q;f(\bx_{0}\cdot\g_{0}\gw))
\left[
\sum_{\g\in\G_{\bx_{0}}(q)}
\gU_{X}(\g\g_{0})
\right]
.
$$
By \eqref{eq:GmodQ}, we have
\be\label{eq:quot}
\G_{\bx_{0}}(q)\bk\G
\cong
\{
(c,d)\in(\Z/q)^{2}:(c,d,q)=1
\}
.
\ee
Apply \eqref{eq:normXq} to the inner brackets, giving 
$$
\cM_{\fq}
=
\cM^{(1)}_{\fq}
+
r^{(1)}(\fq)
,
$$
where
$$
\cM^{(1)}_{\fq}
=
\sum_{q\mid\fq\atop q<Q_{0}}
\gbrho\left(\frac \fq q\right)
\sum_{\gw\in\gW_{Y}}
\sum_{\g_{0}\in\G_{\bx_{0}}(q)\bk\G}
\frac{\Xi(q;f(\bx_{0}\cdot\g_{0}\gw))
}{[\G:\G_{\bx_{0}}(q)]}
\left[
\sum_{\g\in\G}
\gU_{X}(\g)
\right]
,
$$
and 
\beann
|r^{(1)}(\fq)|
&\ll&
\sum_{q\mid\fq\atop q<Q_{0}}
\gbrho\left(\frac \fq q\right)
\sum_{\gw\in\gW_{Y}}
\sum_{\g_{0}\in\G_{\bx_{0}}(q)\bk\G}
|\Xi(q;f(\bx_{0}\cdot\g_{0}\gw))|
\cdot
X^{2\gt}
\\
&\ll&
\fq^{\vep}\
\frac1\fq\
|\gW_{Y}|\
Q_{0}^{3}\
X^{2\gt}
.
\eeann
Here we 
used 
 \eqref{eq:gbDef} and \eqref{eq:XiDef} to
 estimate $|\Xi|\le1$ and $\rho(q)\ll q^{\vep}/q$.
Then \eqref{eq:r1Bnd} follows immediately from \eqref{eq:cAis}.

Returning to $\cM^{(1)}_{\fq}$, we add back in the large factors $q\mid\fq$ and subtract them away, writing
$$
\cM^{(1)}_{\fq}=\cM^{(2)}_{\fq}+r^{(2)}(\fq)
,
$$
say,
where
$$
\cM^{(2)}_{\fq}
=
\sum_{q\mid\fq}
\gbrho\left(\frac \fq q\right)
\sum_{\gw\in\gW_{Y}}
\sum_{\g_{0}\in\G_{\bx_{0}}(q)\bk\G}
\frac{\Xi(q;f(\bx_{0}\cdot\g_{0}\gw))
}{[\G:\G_{\bx_{0}}(q)]}
\left[
\sum_{\g\in\G}
\gU_{X}(\g)
\right]
.
$$
Since $\g_{0}$ ranges over the full quotient, we may drop $\gw$ from $\Xi$, giving
$$
\cM^{(2)}_{\fq}
=
\cX
\rho(\fq)
\sum_{q\mid\fq}
\rho_{1}(q)
=
\cX
\rho(\fq)
\prod_{p\mid\fq}
\left(
1+
\rho_{1}(p)
\right)
.
$$
Here $\cX$ is given by \eqref{eq:cXisis}, and
 $\rho_{1}$ is a multiplicative function defined on primes by
$$
\rho_{1}(p):=
\frac 1 {\rho(p)}
\frac1{[\G:\G_{\bx_{0}}(p)]}
\sum_{\g_{0}\in\G_{\bx_{0}}(p)\bk\G}
\Xi(p;f(\bx_{0}\cdot\g_{0}))
.
$$
A calculation from \eqref{eq:quot} and the definitions shows that 
\be\label{eq:rho1}
\rho_{1}(p)=-{p-1\over(2p-1)(p+1)}
,
\ee
whence
$$
\cM^{(2)}_{\fq}=\gb(\fq)\cX,
$$
with $\gb$ given by \eqref{eq:gbIs}.

It remains to estimate
\beann
r^{(2)}(\fq)
&=&
\cX
\sum_{q\mid\fq\atop q\ge Q_{0}}
\gbrho\left(\frac \fq q\right)
\sum_{\g_{0}\in\G_{\bx_{0}}(q)\bk\G}
\frac{\Xi(q;f(\bx_{0}\cdot\g_{0}\g))
}{[\G:\G_{\bx_{0}}(q)]}
\\
&=&
\cX
\rho(\fq)
\sum_{q\mid\fq\atop q\ge Q_{0}}
\rho_{1}(q)
.
\eeann
It is easy to see that
  $|\rho_{1}(q)|\le1/q$,
and hence we have again that
\beann
|r^{(2)}(\fq)|
&\ll&
\fq^{\vep}
\cX
\frac1\fq
\frac1{Q_{0}}
.
\eeann
Then \eqref{eq:r2Bnd} follows immediately.
\epf

\begin{rmk}\label{rmk:NeedGap}
For our application,  \eqref{eq:r1Bnd} is sufficient as long as 
$$
Q_{0}^{3}
\ll
T^{-\vep}
X^{2(\gd-\gt)}
,
$$
so $Q_{0}$ cannot be too big. But for  \eqref{eq:r2Bnd} to also suffice requires us to take $Q_{0}>T^{\vep_{0}}$. Hence we do need to know that $\G$ has {\it some} spectral gap, but any gap will do; cf. Remark \ref{rmk:prev}. Note also that these error terms pose no restriction on $\cQ$ beyond $\cQ<  T^{C}$.
\end{rmk}

\newpage

\section{Analysis of $\cE$}\label{sec:error}

Recall from \eqref{eq:cEDef} and  \eqref{eq:cAqDec}
that the net error
is given by
$$
\cE
=
\sum_{\fq<\cQ}
\left|
r(\fq)
\right|
,
$$
where
$$
r(\fq)
:=
\sum_{q\mid\fq\atop q\ge Q_{0}}
\sum_{n}a_{T}(n)
\Xi(q;n)\gbrho\left(\frac \fq q\right)
.
$$
The goal of this section is to prove the following
\begin{thm}\label{thm:cEBnd}
Fix any small $\eta>0$, and assume \eqref{eq:YtoT}. Then we have the bound
\be\label{eq:cETotalBnd}
\cE
\ll
T^{\vep}\,
\cX
\cdot
\bigg(
{X^{(1-\gd)}
\over Q_{0}^{\eta/2}}
+
{Y^{2}
\cQ^{6}\over
X^{\gd+5/2}}
\bigg)
.
\ee
\end{thm}

Before proceeding with the proof, we first perform some initial  manipulations to $\cE$. 
Let $\gz(\fq):=\sgn r(\fq)$, and reverse orders of summation, giving 
\beann
\cE
&=&
\sum_{\fq<\cQ}\gz(\fq)\sum_{n}a_{T}(n)\sum_{q\mid\fq\atop q\ge Q_{0}}\Xi(q;n)\gbrho\left(\frac\fq q\right)
\\
&=&
\sum_{Q_{0}\le q<\cQ}
\sum_{n}a_{T}(n)\,
\Xi(q;n)\,
\gz_{1}(q)
,
\eeann
where
$$
\gz_{1}(q)
:=
\sum_{\fq<\cQ\atop \fq\equiv0(q)}
\gz(\fq)
\gbrho\left(\frac\fq q\right)
\ll
T^{\vep}
.
$$
Inserting \eqref{eq:aTDef} gives
$$
\cE=
\sum_{Q_{0}\le q<\cQ}
\sum_{\g\in\G\atop\g=
\bigl( \begin{smallmatrix}
*&*\\ c&d
\end{smallmatrix} \bigr)
}
\gU_{X}(\g)
\sum_{\gw\in\gW_{Y}}
\Xi(q;\ff_{\gw}(c,d))\,
\gz_{1}(q)
.
$$
For $\g=\mattwo**cd\in\G$, the map $\g\mapsto(c,d)$ is $1$-to-$1$ because $\G$ has no parabolic elements. Apply Cauchy-Schwarz in the $\g$ variable, replacing the sum on $\G$ by a smooth sum on $(c,d)\in\Z^{2}$ of norm at most $X$:
$$
\cE^{2}\ll 
X^{2\gd}
\cdot
\sum_{c,d}
\Phi\left(\frac cX\right)
\Phi\left(\frac dX\right)
\left|
\sum_{Q_{0}\le q<\cQ}
\sum_{\gw\in\gW_{Y}}
\Xi(q;\ff_{\gw}(c,d))\,
\gz_{1}(q)
\right|^{2}
.
$$
Here we used \eqref{eq:normX}, and inserted  a fixed function $\Phi$, assumed to be  smooth, real, nonnegative, with 
$\Phi(x)\ge1$ for $x\in[-1,1]$, and   Fourier transform $\hat\Phi$ supported in $[-1,1]$. 
Open the square and reverse orders
$$
\cE^{2}\ll 
X^{2\gd}
\sum_{Q_{0}\le q,q'<\cQ}
T^{\vep}
\left|
\sum_{\gw,\gw'\in\gW_{Y}}
\sum_{c,d}
\Phi\left(\frac cX\right)
\Phi\left(\frac dX\right)
\Xi(q;\ff_{\gw}(c,d))
\Xi(q';\ff_{\gw'}(c,d))
\right|
.
$$
(Recall $\Xi$ is real, so we do not need complex conjugates.) 
Break the sum according to whether $[q,q']\le X$ or $>X$,
writing
\be\label{eq:cE2lg}
\cE^{2}
\ll 
T^{\vep}\,
X^{2\gd}
\cdot
\bigg[
\cE_{\le}
+
\cE_{>}
\bigg]
,
\ee
say.
 The former has complete sums in $c,d$ and the latter is incomplete, so we 
  handle these separately.

\subsection{Range 1: $[q,q']\le X$}\

In this subsection, we will prove the following
\begin{prop}\label{prop:cEle}
Let
$$
\cE_{\le}
:=
\sum_{Q_{0}\le q,q'<\cQ\atop[q,q']\le X}
\left|
\sum_{\gw,\gw'\in\gW_{Y}}
\sum_{c,d}
\Phi\left(\frac cX\right)
\Phi\left(\frac dX\right)
\Xi(q;\ff_{\gw}(c,d))
\Xi(q';\ff_{\gw'}(c,d))
\right|
,
$$
and fix any $\eta>0$ sufficiently small. Then assuming \eqref{eq:YtoT}, we have
\be\label{eq:cEleBnd}
\cE_{\le}
\quad
\ll
\quad
T^{\vep}
X^{2}
|\gW_{Y}|^{2}
{1\over Q_{0}^{\eta}}
.
\ee
\end{prop}

Before the proof, we need two lemmata.
Let
\be\label{eq:S1Def}
S_{1}(q;\gw):=
\frac1{q^{2}}
\sum_{c,d(q)}
\Xi(q;\ff_{\gw}(c,d))
,
\ee
and
\be\label{eq:S2Def}
S_{2}(q;\gw,\gw'):=
\frac1{q^{2}}
\sum_{c,d(q)}
\Xi(q;\ff_{\gw}(c,d))
\Xi(q;\ff_{\gw'}(c,d))
.
\ee


\begin{lem}\label{lem:Dispersion}
If $q>1$, then
$$
S_{1}(q;\gw)=0.
$$
\end{lem}
\pf
Write $\ff_{\gw}(c,d)=Ac^{2}+2Bcd+Ccd,$ and recall that $1=
\left|
\bigl( \begin{smallmatrix}
A&B\\ B&C
\end{smallmatrix} \bigr)
\right|
=AC-B^{2}$.
By multiplicativity, we reduce to the case $q=p$ a prime. Inserting the definitions \eqref{eq:XiDef} and \eqref{eq:gbDef}, we have
\beann
S_{1}(p;\gw)
&=&
\frac1{p^{2}}
\sum_{c,d(p)}
\bo_{\{\ff_{\gw}(c,d)\equiv0(p)\}}
-
{2p-1\over p^{2}}
.
\eeann
We count the number of solutions via exponential sums:
\beann
S_{1}(p;\gw)
&=&
\frac1{p^{2}}
\sum_{c,d(p)}
\frac1p
\sideset{}{'}\sum_{b(p)}e_{p}(b\ff_{\gw}(c,d))
+\frac1p
-
{2p-1\over p^{2}}
,
\eeann
where we separated the term $b=0$. Since $AC-B^{2}=1$, we may assume that $(A,p)=1$. Then completing the square and evaluating Gauss sums gives
\beann
S_{1}(p;\gw)
&=&
\frac1p
\sideset{}{'}\sum_{b(p)}
\frac1{p^{2}}
\sum_{c,d(p)}
e_{p}(bA(c+B\bar Ad)^{2}+b\bar Ad^{2})
-
{p-1\over p^{2}}
\\
&=&
\frac1p
\sideset{}{'}\sum_{b(p)}
\left({aA\over p}\right)p^{-1/2}
\left({a\bar A\over p}\right)p^{-1/2}
-
{p-1\over p^{2}}
=
0
,
\eeann
as claimed.
\epf

\begin{lem}\label{lem:S2}
With $S_{2}$ given by \eqref{eq:S2Def}, we have
$$
|S_{2}(q;\gw,\gw')|
\ll
{q^{\vep}\over q^{2}}
\sum_{q_{1}\mid q}
\left[
\sum_{c,d(q_{1})\atop (c,d,q_{1})=1}
\bo_{\left\{{\ff_{\gw}(c,d)\equiv0(q_{1})\atop\ff_{\gw'}(c,d)\equiv0(q_{1})}\right\}}
\right]
.
$$
\end{lem}
\pf
By multiplicativity we again reduce to the case $q=p$. Then
\beann
S_{2}(p;\gw,\gw')
&=&
\frac1{p^{2}}
\sum_{c,d(p)}
\left(
\bo_{\{\ff_{\gw}(c,d)\equiv0(p)\}}-\gbrho(p)
\right)
\left(
\bo_{\{\ff_{\gw'}(c,d)\equiv0(p)\}}-\gbrho(p)
\right)
\\
&=&
\frac1{p^{2}}
\sum_{c,d(p)\atop(c,d)\not\equiv(0,0)}
\bo_{\left\{{\ff_{\gw}(c,d)\equiv0(p)\atop\ff_{\gw'}(c,d)\equiv0(p)}\right\}}
+
\frac1{p^{2}}
-
\gbrho(p)^{2}
,
\eeann
whence the claim follows.
\epf

We now proceed with the
\pf[Proof of Proposition \ref{prop:cEle}]
Write 
\be\label{eq:tilq}
\tilde q:=(q,q'),\quad q=q_{1}\tilde q 
,\quad q'=q'_{1}\tilde q 
,\quad \bar q:=[q,q']=q_{1}q_{1}'\tilde q 
.
\ee
Applying Poisson summation and splitting the $c,d$ sum by the Chinese Remainder Theorem, we have
\bea
\nonumber
\cE_{\le}
&=&
\sum_{Q_{0}\le q,q'<\cQ\atop \bar q\le X}
\left|
\sum_{\gw,\gw'\in\gW_{Y}}
\sum_{c,d(\bar q)}
\Xi(q;\ff_{\gw}(c,d))
\Xi(q';\ff_{\gw'}(c,d))
\frac{X^{2}}{\bar q^{2}}
\hat\Phi(0)^{2}
\right|
\\
\nonumber
&=&
\hat\Phi(0)^{2}
X^{2}
\sum_{Q_{0}\le q,q'<\cQ\atop \bar q\le X}
\left|
\sum_{\gw,\gw'\in\gW_{Y}}
S_{1}(q_{1};\gw)
S_{1}(q_{1}';\gw')
S_{2}(\tilde q;\gw,\gw')
\right|
,
\eea
with $S_{1}$, $S_{2}$ given by \eqref{eq:S1Def}, \eqref{eq:S2Def}.

Applying Lemma \ref{lem:Dispersion},
 the sum completely vanishes unless $q_{1}=1=q'_{1}$, that is, $q=q'=\tilde q=\bar q$. So we have that
\be\label{eq:cEle2}
\cE_{\le}
=
\hat\Phi(0)^{2}
X^{2}
\sum_{Q_{0}\le q<\cQ}
\left|
\sum_{\gw,\gw'\in\gW_{Y}}
S_{2}(q;\gw,\gw')
\right|
.
\ee

Inserting Lemma \ref{lem:S2} into \eqref{eq:cEle2} gives
$$
\cE_{\le}
\ll
T^{\vep}
X^{2}
\sum_{Q_{0}\le q<\cQ}
\frac1{q^{2}}
\sum_{q_{1}\mid q}
\sum_{c,d(q_{1})\atop(c,d,q_{1})=1}
\sum_{\gw,\gw'\in\gW_{Y}}
\bo_{\left\{{\ff_{\gw}(c,d)\equiv0(q_{1})\atop\ff_{\gw'}(c,d)\equiv0(q_{1})}\right\}}
.
$$

With $c$ and $\gw$ fixed, there are at most two values of $d(\mod q_{1})$ with $\ff_{\gw}(c,d)\equiv0(q_{1})$, and $(c,d,q_{1})=1$. For such $d$, we save $q_{1}^{\eta}$ from the $\gw'$ sum using \eqref{eq:gWYsave}.
In total, we have
\bea
\nonumber
\cE_{\le}
&\ll&
T^{\vep}
X^{2}
\sum_{Q_{0}\le q<\cQ}
\frac1{q^{2}}
\sum_{q_{1}\mid q}
q_{1}
|\gW_{Y}|\,
|\gW_{Y}|q_{1}^{-\eta}
\\
\nonumber 
&\ll&
T^{\vep}
X^{2}
|\gW_{Y}|^{2}
\sum_{Q_{0}\le q<\cQ}
{1\over q^{1+\eta}}
\ll
T^{\vep}
X^{2}
|\gW_{Y}|^{2}
{1\over Q_{0}^{\eta}}
,
\eea
as claimed.
\epf


\subsection{Range 2: $[q,q']>X$}\

Next we study $\cE_{>}$, which recall is given by
$$
\cE_{>}:=
\sum_{Q_{0}\le q,q'<\cQ\atop[q,q']>X}
\left|
\sum_{\gw,\gw'\in\gW_{Y}}
\sum_{c,d}
\Phi\left(\frac cX\right)
\Phi\left(\frac dX\right)
\Xi(q;\ff_{\gw}(c,d))
\Xi(q';\ff_{\gw'}(c,d))
\right|
.
$$ 
Here the $c,d$ sum is incomplete and needs to be completed. 
We will prove the following
\begin{prop}\label{prop:cEg}
We have the estimate
\be
\label{eq:cEgBnd}
\cE_{>}
\ \ll \
\cQ^{\vep}
{\cQ^{12}\over
X^{5}}
|\gW_{Y}|^{2}
Y^{4}
.
\ee
\end{prop}

We again begin with some local lemmata.
Let
$$
S_{4}(q;k,\ell;\gw):=
\frac 1{q^{2}}
\sum_{c,d(q)}
\Xi(q;\ff_{\gw}(c,d))e_{q}(-ck-d\ell)
,
$$
and
$$
S_{5}( q;k,\ell;\gw,\gw')
:=
\frac 1{q^{2}}
\sum_{c,d(q)}
\Xi(q;\ff_{\gw}(c,d))
\Xi(q;\ff_{\gw'}(c,d))
e_{q}(-ck-d\ell)
.
$$

\begin{lem}\label{lem:S4}
We have
$$
|S_{4}(q;k,\ell;\gw)|\le
\twocase{}{
{(
\ff_{\gw}(\ell,-k)
,q)\over q^{2}}
}{if $(k,\ell,q)=1$,}
{0}{if $(k,\ell,q)>1$
.}
$$
\end{lem}
\pf
By multiplicativity, we reduce to the case $q=p$. 
If $k\equiv\ell\equiv0(p)$, then $S_{4}=S_{1}=0$ by Lemma \ref{lem:Dispersion}.

Otherwise, we evaluate
\beann
S_{4}(p;k,\ell;\gw)
&=&
\frac 1{p^{2}}
\sum_{c,d(p)}
\left(
\bo_{\{\ff_{\gw}(c,d)\equiv0(p)\}}
-
\gbrho(p)
\right)
e_{p}(-ck-d\ell)
\\
&=&
\frac 1{p^{2}}
\sum_{c,d(p)}
\frac1p
\sideset{}{'}\sum_{b(p)}
e_{p}(b\ff_{\gw}(c,d))
e_{p}(-ck-d\ell)
,
\eeann
since at least one of $k,\ell$ is non-zero. Assuming without loss of generality that $(A,p)=1$, we have that
\beann
b\ff_{\gw}(c,d))
-ck-d\ell
&\equiv&
bA(c+B\bar Ad-\bar 2\bar b\bar Ak)^{2}
+
b\bar A(d+\bar 2\bar b Bk-\bar 2\bar b A\ell)^{2}
\\
&&
\hskip.5in
-\bar 4\bar b
\ff_{\gw}(\ell,-k)
\qquad(\mod p).
\eeann
Evaluating the Gauss sums in $c$ and $d$ gives
\beann
S_{4}(p;k,\ell;\gw)
&=&
\frac1p
\sideset{}{'}\sum_{b(p)}
\frac1p
e_{p}(-\bar 4\bar b
\ff_{\gw}(\ell,-k)
)
\\
&=&
\twocase{}
{p-1\over p^{2}}
{if $(
\ff_{\gw}(\ell,-k)
,p)=p$,}
{-1\over p^{2}}
{if $(
\ff_{\gw}(\ell,-k)
,p)=1$,}
\eeann
from which the claim follows.
\epf

We treat $S_{5}$ even more trivially than $S_{2}$. 
\begin{lem}\label{lem:S5}
$$
|S_{5}( q;k,\ell;\gw,\gw')|\ll {q^{\vep}\over q}.
$$
\end{lem}
\pf
Again, we consider $q=p$, and estimate trivially
\beann
|S_{5}( p;k,\ell;\gw,\gw')|
&\le&
\frac1{p^{2}}
\sum_{c,d(p)}
\left(
\bo_{\ff_{\gw}(c,d)\equiv0(p)}
+\gbrho(p)
\right)
\le
\frac4p
.
\eeann
The claim then follows.
\epf

\pf[Proof of Proposition \ref{prop:cEg}]
Recalling the notation in \eqref{eq:tilq}, we 
complete the incomplete sums, giving
\be
\label{eq:cEg}
\cE_{>}
=
\sum_{Q_{0}\le q,q'<\cQ\atop\bar q>X}
\left|
\sum_{\gw,\gw'\in\gW_{Y}}
\sum_{0\le k,\ell<\bar q}
\cI(X;k,\ell)S_{3}(q,q';k,\ell;\gw,\gw')
\right|
,
\ee
where
$$
S_{3}(q,q';k,\ell;\gw,\gw')
:=
\frac1{\bar q^{2}}
\sum_{c,d(\bar q)}
\Xi(q;\ff_{\gw}(c,d))
\Xi(q';\ff_{\gw'}(c,d))
e_{\bar q}(-ck- d\ell)
,
$$
and
\bea
\nonumber
\cI(X;k,\ell)
&:=&
\bar q^{2}
\sum_{n,m\in\Z^{2}}
\Phi\left(\frac nX\right)
\Phi\left(\frac mX\right)
e_{\bar q}(nk+m\ell)
\\
\nonumber
&=&
X^{2}
\sum_{n,m}
\hat\Phi\left(X\left(\frac k{\bar q}-n\right)\right)
\hat\Phi\left(X\left(\frac \ell{\bar q}-m\right)\right)
\\
\label{eq:cIbnd}
&\ll&
X^{2}\bo_{\{| k|,|\ell|<\bar q/X\}}
,
\eea
by Poisson summation.

Since $\bar q=q_{1}q_{1}'\tilde q$, the $S_{3}$ sum factors as
$$
S_{3}(q,q';k,\ell;\gw,\gw')=
S_{4}(q_{1};k,\ell;\gw)\cdot
S_{4}(q_{1}';k,\ell;\gw')\cdot
S_{5}(\tilde q;k,\ell;\gw,\gw')
.
$$




Applying \eqref{eq:cIbnd} and Lemmata \ref{lem:S4} and \ref{lem:S5}, 
we obtain
\beann
\cE_{>}
&\ll&
\cQ^{\vep}
X^{2}
\sum_{X<\bar q<\cQ^{2}}\
\sum_{q_{1}q_{1}'\tilde q=\bar q}\
\sum_{\gw,\gw'\in\gW_{Y}}\
\sum_{0\le k,\ell<\bar q/X\atop(k,\ell)\neq(0,0)}
\\
&&
\times
{(
\ff_{\gw}(\ell,-k)
,q_{1})\over q_{1}^{2}}
{(
\ff_{\gw'}(\ell,-k)
,q_{1}')\over (q_{1}')^{2}}
\cdot
\frac1{\tilde q}
.
\eeann 

Observe that since $\bar q>X$, we have
$$
\cQ^{2}>qq'=q_{1}\tilde q q_{1}'\tilde q=\bar q\tilde q>X\tilde q,
$$
and hence $\tilde q<\cQ^{2}/X$. Also, for $\gw\in\gW_{Y}$, we have that 
$$
|A|,|B|,|C|\ll Y^{2},
$$ 
and thus
$$
\ff_{\gw}(\ell,-k)
=
A\ell^{2}-2B\ell k+Ck^{2}
\ll 
\left({\bar q\over X}\right)^{2}Y^{2}
.
$$
Since $(k,\ell)\neq(0,0)$ and $\ff_{\gw}$ is 
definite, we have that $
\ff_{\gw}(\ell,-k)
\neq0$, so can estimate
\bea
\nonumber
\cE_{>}
&\ll&
\cQ^{\vep}
X^{2}
\sum_{\bar q<\cQ^{2}}\
\frac1{\bar q^{2}}
|\gW_{Y}|^{2}
\left({\bar q\over X}\right)^{2}
\left[\left({\bar q\over X}\right)^{2}Y^{2}\right]^{2}
\left({\cQ^{2}\over X}\right)
\\
\nonumber
&\ll&
\cQ^{\vep}
{\cQ^{12}\over
X^{5}}
|\gW_{Y}|^{2}
Y^{4}
,
\eea 
as claimed.
\epf

\begin{rmk}\label{rmk:vep0}
For ease of exposition, we have  not exploited all the cancellation we could  out of the $S_{5}$ sum. In particular, 
we have the bound $|S_{5}|\le |S_{2}|$, so could apply Lemma \ref{lem:S2} and save an extra $\tilde q^{\eta}$ in the $\gw'$ sum above, just as we did in the proof of Proposition \ref{prop:cEle}.
Moreover, we could similarly save $(q_{1}q_{1}')^{\eta'}$ by proving a suitable analogue to \eqref{eq:gWYsave}. This extra savings 
may
get the level just a 
hair
higher, but we will not bother; cf. Remark \ref{rmk:gaMore}.
\end{rmk}


\subsection{Proof of Theorem \ref{thm:cEBnd}}\

Inserting \eqref{eq:cEleBnd} and \eqref{eq:cEgBnd} into \eqref{eq:cE2lg} gives
\beann
\cE^{2}
&\ll &
T^{\vep}\,X^{2\gd}\,|\gW_{Y}|^{2}
\cdot
\bigg[
X^{2(1-\gd)}
{1\over Q_{0}^{\eta}}
+
{\cQ^{12}\over
X^{5+2\gd}}
Y^{4}
\bigg]
,
\eeann
whence the claim follows.

\newpage



\section{Proofs of Theorem \ref{thm:main} and Corollary \ref{cor:main}}\label{sec:proofs}

\subsection{Proof of Theorem \ref{thm:main}}\

Let $\cA$ be the sequence given by \eqref{eq:aTDef}. For some small $\eta>0$ to be chosen later, assume \eqref{eq:YtoT}. Collecting \eqref{eq:cAqDec}, \eqref{eq:gbIs}--\eqref{eq:r2Bnd}, and \eqref{eq:cETotalBnd}, we see that
$$
|\cA_{\fq}|=\gb(\fq)\cX+r^{(3)}(\fq),
$$
with \eqref{eq:gbBnd} satisfied by classical methods, and 
\beann
\sum_{\fq<\cQ}|r^{(3)}(\fq)|
\ll
T^{\vep}
\cX
\left[
{Q_{0}^{3}\over T^{2(\gd-\gt)}}
+
\frac1{Q_{0}}
+
{X^{1-\gd}\over Q_{0}^{\eta/2}}
+
{Y^{2}\cQ^{6}\over X^{\gd+5/2}}
\right]
.
\eeann

Write
$$
X=T^{x}
,
Y=T^{y}
,
Q_{0}=T^{2\ga_{0}}
,
\text{ and }
\cQ=T^{2\ga}
,
$$
with $x+y=1$. For \eqref{eq:rqlcA} and \eqref{eq:YtoT} to hold, we need the following inequalities:
\bea
\label{eq:6281}
6\ga_{0}
&<&
2(\gd-\gt)x
\\
\label{eq:6282}
\ga_{0}
&>&
0
\\
\label{eq:6283}
(1-\gd)x
&<&
\ga_{0}\eta
\\
\label{eq:6284}
2y+12\ga
&<&
(\gd+5/2)x
\\
\label{eq:6285}
y
&<&
3\ga\eta(\cC+1)/(\gd-\gt)
.
\eea
\begin{rmk}\label{rmk:NoSpecGap}
With $\gd$ very near $1$ and $y$ very near $0$ (hence $x$ near $1$), it is clear that \eqref{eq:6284} will determine the exponent $\ga$. 
Notice that this condition does not depend on the spectral gap $\gT$, cf. Remark \ref{rmk:prev}.
Heuristically, we should take $y$ tiny, and $\eta$ even much smaller so that \eqref{eq:6285} holds; then \eqref{eq:6283} determines how close $\gd$ must be to $1$, and \eqref{eq:6281}, \eqref{eq:6282} pose no serious restriction. Let us make this precise.
\end{rmk}

Now we fix $\vep>0$ and set 
$$
\ga=\frac7{24}-\vep,
$$ 
as required for Theorem \ref{thm:main}. We will take $\gd$ very near $1$, so may already assume that $\gd>14/15$. Then using Gamburd's gap $\gt=5/6$ (see Theorem \ref{thm:spec}), we have
$$
\frac1{10}
<
\gd-\gt
<\frac16.
$$
Since $X$ is treated as the large variable, assume $x>1/2$. Then we can set 
$$
\ga_{0}=\frac1{100}
,
$$ 
whence \eqref{eq:6281} and \eqref{eq:6282} are easily satisfied. 
Using $x=1-y$, rewrite \eqref{eq:6284} as
$$
\frac7{24}-\vep
=
\ga
<
\frac7{24}
-{1\over12}(1-\gd)
-
\left(
\frac7{24}
+
\frac1{12}(1-\gd)
+\frac16
\right)
y
.
$$
Hence it will suffice to make $1-\gd<\vep$ and $y<\vep$, say. We will soon impose much more stringent restrictions on $\gd$, so focus on $y$. 

Let us set
$$
y=\foh\vep,
$$
say.
Then
 we can take
$$
\eta=10^{-13}\vep
,
$$
so that \eqref{eq:6285} is satisfied  using \eqref{eq:cCis}. 

It only remains to ensure that \eqref{eq:6283} holds. The variables $\ga_{0}$, $\eta$,  and $x=1-y$ are now all determined, so this is a restriction on $\gd$. It is easy to see that
$$
1-\gd<10^{-\gdOvep}\vep
$$
suffices, so we set
$$
\gd_{0}(\vep)=1-10^{-\gdOvep}\vep
,
$$
see Remark \ref{rmk:delVals}.
This completes the proof of Theorem \ref{thm:main}.


\subsection{Proof of Corollary \ref{cor:main}}\

The corollary follows easily from the theorem. 
Since $1/\ga\in(3,4)$, we are in position to capture $R=4$ almost primes. 

By the best available linear weighted almost-prime 
sieve due to Greaves \cite[(1.4)]{Greaves1986}
, we 
can produce $R=4$ almost primes as long as
$$
\ga>\frac1{4-0.103974}
\approx
0.256672
.
$$
Setting $\ga=\frac7{24}-\vep$, we may take $\vep$ as large as
$$
\vep=\frac3{100}.
$$
Hence 
$$
\gd_{0}=1-10^{-\gdO}
$$
suffices, as claimed in Remark \ref{rmk:delVals}. This completes the proof.

\newpage

\bibliographystyle{alpha}

\bibliography{AKbibliog}

\end{document}